\documentclass[11pt]{article}
\usepackage{amssymb}
\newcommand{\be}{\begin{equation}}
\newcommand{\ee}{\end{equation}}
\newcommand{\bea}{\begin{eqnarray}}
\newcommand{\eea}{\end{eqnarray}}
\newcommand{\binom}[2]{{#1 \choose #2}}

\begin{document}


\title{On the Sylvester waves in partition function}
\author{Boris Y. Rubinstein
\\Stowers Institute for Medical Research
\\1000 50th St., Kansas City, MO 64110, U.S.A.}
\date{\today}

\maketitle
\begin{abstract}
Sylvester showed that the partition
function 
can be written as a sum of the polynomial term 
and quasiperiodic components called the Sylvester waves.
Recently an explicit 
expression of the Sylvester wave 
as a finite sum over the 
Bernoulli polynomials of higher order with periodic coefficients was found.
This expression can be also written as the weighted sum of the 
polynomial terms with shifted arguments and this manuscript presents a formal proof for
validity of such representation.
\end{abstract}


\vskip0.4cm
The problem of partitions of positive integers has long history started from
the work of Euler who laid a foundation of the theory of partitions
\cite{GAndrews}, introducing the idea of generating functions.
A. Cayley \cite{Cayley1855} and J.J. Sylvester \cite{Sylv1} provided a new insight and made a
remarkable progress in this field. Sylvester found \cite{Sylv1,Sylv2} the
procedure for computation of a {\it restricted} partition function,
while Cayley described the symmetry properties of these functions. The restricted
partition function $W(s,{\bf d}^m) \equiv W(s,\{d_1,d_2,\ldots,d_m\})$ is a
number of partitions of $s$ into positive integers  $\{d_1,d_2,\ldots,d_m\}$.
Sylvester showed that the partition
function can be split into polynomial and quasiperiodic parts and presented as a sum of components
called the {\em Sylvester waves}
\be
W(s,{\bf d}^m) = \sum_{j=1} W_j(s,{\bf d}^m)\;,
\label{SylvWavesExpand}
\ee
where summation runs over all divisors
of the elements in the set ${\bf d}^m$.
The wave $W_j(s,{\bf d}^m)$ is a quasipolynomial in $s$
closely related to prime roots $\rho_j$ of unity.
Namely, Sylvester showed \cite{Sylv2} that the wave
$W_j(s,{\bf d}^m)$ is a coefficient of
${t}^{-1}$ in the series expansion in ascending powers of $t$ of
\be
F_j(s,t)=\sum_{\rho_j} \frac{\rho_j^{-s} e^{st}}{\prod_{k=1}^{m}
        \left(1-\rho_j^{d_k} e^{-d_k t}\right)}\;,
\quad
\rho_j=\exp(2\pi i n/j),
\label{generatorWj}
\ee
where the summation is made over all prime roots of unity
$\rho_j$ for $n$ relatively prime to $j$
(including unity) and smaller than $j$.
The relation (\ref{generatorWj}) is just a recipe for calculation of the partition function and it
does not provide an explicit expression for $W_j(s,{\bf d}^m)$.

Using the Sylvester recipe we found in \cite{Rub04} a formula for the
Sylvester wave $W_j(s,{\bf d}^m)$ as a finite sum of the Bernoulli
polynomials of higher
order \cite{bat53,Norlund1924} multiplied by a periodic function of integer 
period $j$. 
The polynomial part $W_1(s,{\bf d}^m)$ of the partition function reads
\begin{equation}
W_1(s,{\bf d}^m) =
\frac{1}{(m-1)!\;\pi_m}
B_{m-1}^{(m)}(s + s_m, {\bf d}^m)\;,
\quad
s_m = \sum_{i=1}^m d_i, \ \ \pi_m =  \prod_{i=1}^m d_i\;,
\label{W_1}
\end{equation}
with the Bernoulli polynomials of higher 
order are defined by the generating function \cite{Norlund1924}:
\be
\frac{e^{st} t^m \pi_m}{\prod_{i=1}^m (e^{d_it}-1)} =
\sum_{n=0}^{\infty} B^{(m)}_n(s,{\bf d}^m)
\frac{t^{n}}{n!}\;,
\label{genfuncBernoulli0}
\ee
and some of their properties are described in the Appendix.
The Sylvester wave $W_j(s,{\bf d}^m)$ for $j>1$ can be written as
\bea
W_j(s,{\bf d}^m)  &=& 
\frac{j^{k_j-m}}{(m-1)! \; \pi_{m}}
\sum_{n=0}^{k_j-1} C^{m-1}_n B_n(s+s_m,{\bf d}^{k_j})
\nonumber \\
&\times&
\sum_{{\bf r}=0}^{j-1}
B^{(m)}_{m-n-1}({\bf r}\cdot{\bf d}^{m-k_j},j{\bf d}^{m-k_j})
\Psi_j(s+s_m+{\bf r}\cdot{\bf d}^{m-k_j}),
\label{WjOrig}
\eea
where we introduce two subsets -- 
the subset ${\bf d}^{k_j}=\{d_{1},d_{2},\ldots,d_{k_j}\}$ of the generators $d_i$ divisible by $j$
and the subset
${\bf d}^{m-k_j}=\{d_{k_j+1},d_{k_j+2},\ldots,d_{m}\}$ of those nondivisible by $j$.
The $(m-k_j)$-dimensional vector ${\bf r}$
has the form ${\bf r} = \{r_{k_j+1},r_{k_j+2},\ldots,r_{m}\}$
with $0 \le r_i \le j-1$.
The prime circulator $\Psi_j$ in (\ref{WjOrig}) is the $j$-periodic function
introduced in \cite{Cayley1855}
that can be written as a sum of the simple periodic functions
\be
\Psi_j(s) = 
\sum_{\rho_j} \rho_j^s  = \sum_{n} \Psi_{j,n}(s),
\quad
\Psi_{j,n}(s) = \exp(2\pi i n s/j),
\quad
\sum_{k=0}^{j-1} \Psi_{j,n}(k) = 0,
\label{gencirc}
\ee
for all $1 \le n \le j-1$ relatively prime to $j$.

Considering  (\ref{WjOrig}) it was suggested in \cite{Rub04} to extend the outer summation 
to $m-1$ that with the help of (\ref{A2}) enables to present $W_j(s,{\bf d}^m)$ 
as follows
\bea
W_j(s,{\bf d}^m)  =
\frac{j^{k_j-m}}{(m-1)! \; \pi_{m}}
\sum_{{\bf r}=0}^{j-1}
B^{(m)}_{m-1}(s+s_m+ {\bf r}\cdot{\bf d}^{m-k_j},{\bf d}^{m}_j)
\Psi_j(s+s_m+{\bf r}\cdot{\bf d}^{m-k_j}),
\label{WjFin} 
\eea
where we introduce the modified set 
${\bf d}^{m}_j = {\bf d}^{k_j} \cup j{\bf d}^{m-k_j}$.
Comparing (\ref{WjFin}) to (\ref{W_1}) we observe that the 
Sylvester wave $W_j(s,{\bf d}^m)$ can be viewed as a weighted sum of the 
polynomial part with shifted argument with the weights given by 
the prime circulator $\Psi_j$.
This result is valid in the assumption that the added terms corresponding to 
the values of $k_j \le n \le m-1$ in  (\ref{WjOrig}) sum up to zero. 
It was positively tested on multiple examples but the proof was not given.
The aim of this manuscript is to present the formal proof of added terms vanishing
which indicates that the form (\ref{WjFin}) is indeed correct.



Given positive integers ${\mu},{\nu},j$ define the integer vectors ${\bf d}=\{d_1,d_2,\ldots,d_{\mu}\}$,  
${\bf r}=\{r_1,r_2,\ldots,r_{\mu}\}$ with positive $d_i > 0$ and $0 \le r_i \le j-1$. 
Introduce an arbitrary
vector ${\bf e}=\{e_1,e_2,\ldots,e_q\}$ and the $j$-periodic function $\Psi_{j}(s)$
to define a quantity 
\be
\sigma_{\nu}(s,t,{\bf d},{\bf e}) = \sum_{{\bf r}=0}^{j-1}  
B^{(q)}_{\nu}(t+{\bf r}\cdot{\bf d},{\bf e})  \Psi_{j}(s+{\bf r}\cdot{\bf d}),
\quad
\sum_{k=0}^{j-1}\Psi_{j}(k) = 0,
\label{3_a}
\ee
where $B^{(q)}_{\nu}(s,{\bf e})$ is the Bernoulli polynomial of the 
higher order satisfying the relations (\ref{A1},\ref{A2}).
Use (\ref{gencirc}) to rewrite $\Psi_j(s+{\bf r}\cdot{\bf d})$
as the sum 
\be
\Psi_j(s+{\bf r}\cdot{\bf d})
=\sum_{n} \Psi_{j,n}(s+{\bf r}\cdot{\bf d}) = 
\sum_{n}\Psi_{j,n}(s)  \prod_{i=1}^{m} \Psi_{j,n}(r_i d_i),
\quad\quad
\sum_{k=0}^{j-1}\Psi_{j,n}(k) = 0.
\label{3_aa}
\ee
The polynomial $B^{(q)}_{\nu}(t+{\bf r}\cdot{\bf d},{\bf e})$ satisfies the relations
\be
B_{\nu}^{(q)}(t+{\bf r}\cdot{\bf d},{\bf e}) 
=\sum_{p=0}^{\nu} C_p^{\nu}
B_{\nu-p}^{(q)}(t,{\bf e})
B_p({\bf r}\cdot{\bf d}),
\quad
B_{p}({\bf r}\cdot{\bf d}) 
=\sum_{k=0}^{p} C_k^{p}
({\bf r}\cdot{\bf d})^{k}B_{p-k}\;,
\label{3_b}
\ee
where 
the factor $({\bf r}\cdot{\bf d})^{k}$ can be written
as 
\be
({\bf r}\cdot{\bf d})^{k} = 
\left(
\sum_{i=1}^{\mu} r_i d_i
\right)^{k}
=
\sum_{\bf l}  C^{k}_{\bf l}\prod_{i=1}^{\mu} (r_i d_i)^{l_i},
\quad
C^{k}_{\bf l} = \frac{k!}{l_1!l_2!\ldots l_{\mu}!},
\quad
\sum_{i=1}^{\mu} l_i = k,
\label{3_c}
\ee
and $C^{k}_{\bf l}$ denotes the multinomial coefficient.
Substitute (\ref{3_c}) into (\ref{3_b}) to obtain
\be
B_{n}^{(q)}(t+{\bf r}\cdot{\bf d},{\bf e}) 
=\sum_{p=0}^{n} \sum_{k=0}^{p}  C_p^{\nu}C_k^{p}
B_{n-p}^{(q)}(t,{\bf e})B_{p-k}
\sum_{\bf l}  C^{k}_{\bf l}\prod_{i=1}^{\mu} (r_i d_i)^{l_i}\;.
\label{3_d}
\ee
Use it together with  (\ref{3_aa}) in the definition (\ref{3_a}) 
\be
\sigma_{\nu}(s,t,{\bf d},{\bf e}) = 
\sum_{p=0}^{\nu} \sum_{k=0}^{p} C_p^{\nu}C_k^{p}
B_{\nu-p}^{(q)}(t,{\bf e})B_{p-k}
\sum_{n}\Psi_{j,n}(s) 
\sum_{\bf l}  C^{k}_{\bf l}
\sum_{{\bf r}=0}^{j-1}\prod_{i=1}^{\mu}  (r_i d_i)^{l_i} 
\Psi_{j,n}(r_i d_i).
\label{3_e}
\ee
Note that for $k < \mu$  {\it all} terms in the sum over ${\bf l}$ 
have at least one $l_i$ equal to zero and the corresponding
factor $(r_i d_i)^{l_i}$ is absent.
In this case the sum over ${\bf r}$
$$
\sum_{{\bf r}=0}^{j-1} 
\prod_{i=1}^{\mu} (r_i d_i)^{l_i} \Psi_{j,n}(r_i d_i) = 0,
$$
due to (\ref{3_aa}) and the corresponding contribution to $\sigma_{\nu}(s,t,{\bf d},{\bf e})$ vanishes.
It follows from (\ref{3_b}) that the maximal value of $k$ is equal to $\nu$ and thus
the function $\sigma_{\nu}(s,t,{\bf d},{\bf e})$ satisfies
\be
\sigma_{\nu}(s,t,{\bf d},{\bf e}) = 0, 
\qquad
0 \le \nu \le \mu-1.
\label{3_f}
\ee
Consider the inner sum in (\ref{WjOrig}) that can be written as 
$$
\sum_{{\bf r}=0}^{j-1}
B^{(m)}_{m-n-1}({\bf r}\cdot{\bf d}^{m-k_j},j{\bf d}^{m-k_j})
\Psi_j(s+s_m+{\bf r}\cdot{\bf d}^{m-k_j}) 
= \sigma_{m-n-1}(s+s_m,0,{\bf d}^{m-k_j},j{\bf d}^{m-k_j}),
$$
so that $\nu = m-n-1$ and $\mu = m-k_j$. 
Use (\ref{3_f}) to observe that for 
$$
0 \le m-n-1 \le m-k_j-1 
\quad
\Rightarrow
\quad
k_j \le n \le m-1,
$$
the sum $ \sigma_{m-n-1}(s+s_m,0,{\bf d}^{m-k_j},j{\bf d}^{m-k_j})$ vanishes.
This means that the addition to (\ref{WjOrig}) the terms required to produce (\ref{WjFin}) 
does not change the value of the Sylvester wave $W_j(s,{\bf d}^m)$.

Replacing in the r.h.s. of (\ref{WjFin}) all the factors $\Psi_j(s+s_m+{\bf r}\cdot{\bf d}^{m-k_j})$ by unity and
using the multiplication theorem (\ref{A3})
we find
\bea
\frac{j^{k_j-m}}{(m-1)! \; \pi_{m}}
\sum_{{\bf r}=0}^{j-1}
B^{(m)}_{m-1}(s+s_m+ {\bf r}\cdot{\bf d}^{m-k_j},{\bf d}^{m}_j) = 
\frac{1}{(m-1)! \; \pi_{m}} 
B^{(m)}_{m-1}(s+s_m,{\bf d}^{m}) = W_1(s,{\bf d}^m).
\label{unit_weight}
\eea
This relation means that the Sylvester wave $W_j(s,{\bf d}^m)$
represents the splitting of the polynomial part $W_1(s,{\bf d}^m)$ of the partition function
into a multiple sum of $j$-periodic contributions.


\appendix

\section*{Appendix \label{appendix1}}
\renewcommand{\theequation}{A\arabic{equation}}
\setcounter{equation}{0}


The symbolic technique for manipulating sums with binomial coefficients by
expanding polynomials and then replacing powers by subscripts was developed in
nineteenth century by Blissard,
it is known as the umbral
calculus \cite{Roman1978}.
An example of this notation is also found in \cite{bat53} in
section devoted to the Bernoulli polynomials $B_k(x)$.

The well-known formulas for the ordinary Bernoulli polynomials
can be  written symbolically 
\be
B_n(x) = \sum_{k=0}^{n} C^n_k
B_k x^{n-k} \equiv (B+x)^n,\ \
B_n(x+y) = \sum_{k=0}^{n} C^n_k
B_k(x) y^{n-k} \equiv (B(x)+y)^n,
\label{A1}
\ee
where
after the expansion the exponents of $B$ and $B(x)$ are converted into the
orders of the Bernoulli number and Bernoulli polynomial, respectively:
$B^k \Rightarrow B_k, B^k(x) \Rightarrow B_k(x)$.

N\"orlund \cite{Norlund1924} introduced the Bernoulli polynomials of higher order 
having an umbral representation
$$
B_{n}^{(m)}(x,{\bf d}^m) =
\left(
d_m B + B^{(m-1)}(x,{\bf d}^{m-1})
\right)^n,
$$
that recursively reduces to a more symmetric form
\be
B_{n}^{(m)}(x,{\bf d}^m) =
\left(
x + d_1  B +
d_2   B + \ldots + d_m  B
\right)^n =
\left(
x + \sum_{i=1}^m d_i   B
\right)^n.
\label{A2}
\ee
These polynomials satisfy  \cite{Norlund1924} the multiplication theorem
\be
\sum_{r_1=0}^{m_1-1}
\ldots
\sum_{r_p=0}^{m_p-1}
B^{(n)}_{k}(s+\sum_{i=1}^p \frac{r_id_i}{m_i},
\{d_1,\ldots,d_n\})
=
\prod_{i=1}^p m_i \cdot
B^{(n)}_{k}(s,\{\frac{d_1}{m_1},\ldots,\frac{d_p}{m_p},d_{p+1},\ldots,d_n\}).
\label{A3}
\ee



\begin{thebibliography}{99}
\bibitem{GAndrews}
G. E. Andrews, {\it The Theory of Partitions}, Encyclopedia of Mathematics and
its Applications, V.2, 
Addison--Wesley, 1976.
\bibitem{bat53}
H. Bateman and A. Erdel\'yi, {\it Higher Transcendental Functions}, V.1, McGraw-Hill Book Co., NY, 1953.
\bibitem{Cayley1855}
A. Cayley, {\it Researches on the partitions of numbers},
Phil. Trans. Roy. Soc. {\bf CXLV} (1855), 127-140; 
{\it Coll. Math. Papers}, Cambridge Univ. Press, V.II, 235-249, 1889.
\bibitem{Norlund1924}
N.E. N\"orlund, {\it Vorlesungen \"Uber Differenzenrechnung},
Verlag von Julius Springer, Berlin, 1924.
\bibitem{Roman1978}
S. Roman and G.-C. Rota, {\it The umbral calculus}, 
Adv. Math. {\bf 27} (1978), 95-188.
\bibitem{Rub04} B.Y. Rubinstein,
{\it Expression for restricted partition function through
Bernoulli polynomials}, 
Ramanujan Journal {\bf 15} (2008), 177-185.
\bibitem{Sylv1}
J.J. Sylvester, {\it On the partition of numbers},
Quarterly J. Math. {\bf I} (1857), 141-152; 
{\it Coll. Math. Papers}, Cambridge Univ. Press, V.II, 90-99, 1908.
\bibitem{Sylv2}
J. J. Sylvester, {\it On subinvariants, that is, semi-invariants to binary quantics
of an unlimited order},
American Journal of Mathematics {\bf V} (1882), 79-136;
{\it Coll. Math. Papers}, Cambridge Univ. Press, V.III, 568-622, 1909.
\end{thebibliography}
\end{document}